\documentclass[%
 reprint,
 amsmath,amssymb,
 aps,
]{revtex4-1}

\usepackage{graphicx}
\usepackage{dcolumn}
\usepackage{bm}

\usepackage{graphicx,color,mathtools,amsmath,amssymb,float,xcolor,bbold,subfigure}

\def\ul{\underline}

\begin{document}


\title{Using Noisy or Incomplete Data to Discover Models of Spatiotemporal Dynamics}

\author{Patrick A.K. Reinbold}
\author{Daniel R. Gurevich}
 \author{Roman O. Grigoriev}
\affiliation{
School of Physics, Georgia Institute of Technology, \\
 Atlanta, Georgia 30332-0430, USA
}

\date{\today}

\begin{abstract}
Sparse regression has recently emerged as an attractive approach for discovering models of spatiotemporally complex dynamics directly from data. In many instances, such models are in the form of nonlinear partial differential equations (PDEs); hence sparse regression typically requires evaluation of various partial derivatives. However, accurate evaluation of derivatives, especially of high order, is infeasible when the data are noisy, which has a dramatic negative effect on the result of regression. 
We present a novel and rather general approach that addresses this difficulty by using a weak formulation of the problem.
For instance, it allows accurate reconstruction of PDEs involving high-order derivatives, such as the Kuramoto-Sivashinsky equation, from data with a considerable amount of noise. 
The flexibility of our approach also allows reconstruction of PDE models that involve latent variables which cannot be measured directly with acceptable accuracy. This is illustrated by reconstructing a model for a weakly turbulent flow in a thin fluid layer, where neither the forcing nor the pressure field is known. 
\end{abstract}

\pacs{Valid PACS appear here}
\maketitle


\section*{\label{sec:intro}Introduction}

Macroscopic description of numerous physical, chemical, and biological systems typically involves one or several partial differential equations (PDEs). In some instances, these PDEs represent a physical conservation law, and in others, the PDEs are obtained by homogenization of an underlying microscopic description. The Navier-Stokes equation governing fluid flow and the diffusion equation governing heat or mass flux are examples that incorporate both approaches. Despite the differences in their origin, one thing remained constant for several centuries: PDE models were mainly derived from first principles. Their coefficients typically involve either fundamental physical constants, such as the Planck constant in the Schr\"odinger equation, or properties of the system, such as fluid viscosity or thermal conductivity, that can be either computed or measured independently.

An alternative approach -- data-driven discovery of mathematical models, where both the form of the model and the values of the coefficients are determined based solely on the available data -- has emerged relatively recently \cite{crutchfield1987,bongard2007,yao2007,chou2009,brunton2016}. 
In particular, sparse symbolic regression \cite{xu_2008,rudy2017,li_2019} has been applied successfully to identifying PDE models from data with minimal noise (i.e., standard deviation of 1\% or less of the data range).
Unfortunately, since existing approaches based on sparse regression rely on explicit evaluation of various candidate terms using local data, they all experience serious difficulties in the presence of higher noise levels characteristic of typical experimental measurements and generally fail to reconstruct PDE models involving higher order derivatives. 

Another limitation of existing approaches is that they require that all the variables present in the model be either directly observable or local functions of the directly observable data. For instance, using direct measurements of the fluid velocity ${\bf u}$, it is possible to reconstruct the vorticity equation \cite{rudy2017}, which involves ${\bf u}$ and the vorticity $\omega=\nabla\times{\bf u}$, but not the Navier-Stokes equation, which involves ${\bf u}$ and a latent variable -- pressure. 
A recently introduced extension of the sparse regression method circumvents the latter limitation at the expense of raising the order of all of the derivatives \cite{reinbold_2019}. 

An alternative approach that treats time evolution as a Gaussian process \cite{raissi2018siam} was shown to be capable of reconstructing the coefficients in the 2D Navier-Stokes equation without using the pressure field \cite{raissi_2018}. 
However, this approach assumes the model to be known {\it a priori} and exhibits noise sensitivity similar to that of sparse regression-based approaches.
To the best of the authors' knowledge, no method currently exists that can robustly reconstruct PDEs involving latent variables (i.e., variables that cannot be measured) and/or derivatives of a high order using data with high levels of noise, which significantly limits the practical utility of data-driven approach to model discovery.

The present article removes the major roadblock for the data-driven approach in reconstructing PDE-based mathematical models by introducing a weak formulation of the sparse regression method, which addresses both of the limitations mentioned previously. In the following, we introduce the mathematical foundations of our approach and illustrate it using three representative examples: the Kuramoto-Sivashinsky equation, a quasi-two-dimensional fluid flow, and the $\lambda-\omega$ reaction-diffusion system.

\section*{\label{sec:moddisc}Data-Driven Model Discovery}
Models of continuous spatially distributed systems tend to have the form of a PDE
\begin{align}
\sum_{n=0}^Nc_n{\bf f}_n({\bf u},\partial_t{\bf u},\partial_t^2{\bf u},\nabla{\bf u},\nabla^2{\bf u},\cdots)=0,
\label{eq:nl1}
\end{align}
where each of the terms depends on the system state ${\bf u}$ and its spatial and temporal derivatives of various orders and $c_n$ are coefficients assumed to be constant in this study (an extension to coefficients depending on spatial and/or temporal coordinates is straightforward \cite{xu_2008,li_2019}). 
Symmetry and physical constraints can be used to narrow down the functional form of the terms that can appear in the model \cite{reinbold_2019}, and sparse regression can be used to discard unnecessary terms and determine a parsimonious form of the model and the values of the corresponding coefficients $c_n$.

We will illustrate the procedure using examples that involve a single term containing a temporal derivative of the state ${\bf u}$. The corresponding coefficient can be set to unity without loss of generality. Separating this term on the left-hand-side, we can rewrite \eqref{eq:nl1} as
\begin{align}
\partial_t^k\hat{D}{\bf u}=\sum_{n=1}^Nc_n{\bf f}_n({\bf u},\nabla{\bf u},\nabla^2{\bf u},\cdots),
\label{eq:nl2}
\end{align}
where $\hat{D}$ is typically either an identity or a linear operator involving only spatial derivatives and $k$ is the order of the temporal derivative. For instance, $k=1$ and $\hat{D}=\mathbb{1}$ for the Navier-Stokes equation, $k=2$ and $\hat{D}=\mathbb{1}$ for the wave equation, $k=1$ and $\hat{D}=\nabla^2$ for the Orr-Sommerfeld equation, etc.

To convert this to a linear algebra problem amenable to sparse regression, let us multiply the differential equation \eqref{eq:nl2} by a weight ${\bf w}$ and integrate the result over a spatiotemporal domain $\Omega_k$, then repeat this procedure for $K$ different choices of $\Omega_k$. This will yield a system 
\begin{align}
{\bf q}_0 = \sum_{n=1}^N c_n{\bf q}_n = Q{\bf c},
\label{eq:symreg}
\end{align}
where $Q=[{\bf q}_1 \cdots {\bf q}_N]$ is the ``library'' and the ``library terms'' ${\bf q}_n\in\mathbb{R}^K$ are column vectors corresponding to different terms $f_n$ in \eqref{eq:nl2} with entries that correspond to a particular choice of ${\bf w}$ and $\Omega_k$, e.g.,
\begin{align}
q_n^k = \int_{\Omega_k}{\bf w}\cdot{\bf f}_n\,d\Omega.
\label{eq:q}
\end{align}
The key advantage of this formulation compared to the local approach investigated previously \cite{xu_2008,rudy2017,li_2019,reinbold_2019} is that, by performing integration by parts, the action of derivatives can be transferred from the noisy data ${\bf u}$ onto the smooth weight ${\bf w}$, dramatically decreasing the effect of noise on terms involving high-order derivatives. Furthermore, the weight function can be chosen in such a way that the terms involving latent variables are eliminated, yielding a problem that can be solved using standard techniques.

A parsimonious model can finally be determined by choosing $K\ge N$ and using an iterative sparse regression algorithm such as SINDy \cite{brunton2016}. Each iteration involves computing the solution 
\begin{align}\label{eq:tildec}
\tilde{\bf c}= Q^+{\bf 	q}_0,
\end{align}
which minimizes the residual of the linear system defined by \eqref{eq:symreg}, where $Q^+$ denotes the pseudo-inverse of $Q$.
This is followed by a thresholding procedure to remove dynamically irrelevant terms with $\|\tilde{c}_n{\bf q}_n\|<\gamma\|{\bf q}_0\|$ for sufficiently small $\gamma$ (we choose $\gamma=0.05$).
To validate the results of regression, we use an ensemble of $M$ cases with different random distributions of the $K$ integration domains $\Omega_k$ relative to the spatiotemporal domain on which the data are available (we use $M=30$ and $K=100$).

Our approach is illustrated below using several examples that highlight different aspects of the problem. 
In the first two, we will assume that the form of the model is known, so that only the coefficients have to be determined. 
The last example illustrates how a parsimonious model can be identified via symbolic regression using a large library of candidate terms.
In each case, we generate the surrogate data using the reference nonlinear PDE, add noise with standard deviation $\sigma$ to this data, evaluate the integrals using the composite trapezoidal rule, and then solve the sparse regression problem to reconstruct the reference PDE.
Note that, in all cases, the range of the data is $O(1)$, so that $\sigma=1$ corresponds to 100\% noise. Numerical codes used to generate the datasets and the MATLAB codes used to identify the governing equations using these datasets can be found in the GitHub repository: \url{https://github.com/pakreinbold/PDE_Discovery_Weak_Formulation}.

\section*{\label{sec:KS} High order derivatives}

The Kuramoto-Sivashinsky equation
\begin{align}
\partial_tu + u\partial_xu + \partial^2_xu + \partial^4_xu = 0,
\label{eq:KS}
\end{align}
describes the chaotic dynamics of laminar flame fronts \cite{sivashinsky_1977}, reaction-diffusion systems \cite{kuramoto_1976}, and coating flows \cite{sivashinsky_1980flow}.
This is a notable example of a nonlinear PDE that involves high-order partial derivatives, which has made it difficult to accurately reconstruct from noisy data.
Rearranging this PDE into the form of Eq. \eqref{eq:nl2}, we find  $c_1=c_2=c_3=-1$.

Since the Kuramoto-Siva\-shinsky equation involves a scalar variable $u$, it can be converted to weak form by integrating its product with a scalar weight $w$ over a set of different integration domains
\begin{align}
\Omega_k=\{(x,t):|x-x_k|\le H_x,|t-t_k|\le H_t\}
\end{align}
centered around randomly chosen points $(x_k,t_k)$.
This yields the system \eqref{eq:symreg} with library terms whose elements are given by
\begin{align}
q_0^k  & = \int_{\Omega_k} w \partial_t u\,d\Omega, & q_1^k & = \int_{\Omega_k} wu\partial_xu\,d\Omega, \nonumber\\
q_2^k & = \int_{\Omega_k} w \partial^2_xu\,d\Omega, & q_3^k & = \int_{\Omega_k} w \partial^4_xu\,d\Omega.
\end{align}

Integration by parts can be used to move all derivatives from the noisy field $u$ onto a smooth noiseless $w$, yielding
\begin{align}
q^k_0  & = -\int_{\Omega_k} u \partial_t w \,d\Omega, & q^k_1 & = -\int_{\Omega_k} \frac{1}{2}u^2\partial_xA\,d\Omega, \nonumber\\ 
q^k_2 & = \int_{\Omega_k} u \partial^2_xw \,d\Omega, & q^k_3 & = \int_{\Omega_k} u \partial^4_xw\,d\Omega,
\label{eq:KSlib}
\end{align}
provided $w$ satisfies the conditions required for the boundary terms to vanish. Specifically, $w$ (and its derivatives up to third order in space) should vanish along the boundary $\partial\Omega_k$. To satisfy these boundary conditions, we chose 
\begin{align}
w = (\ul{x}^2-1)^p(\ul{t}^2-1)^q,
\end{align}
where $p\ge 4$, $q\ge 1$ are integers and the underbar denotes nondimensionalized variables $\ul{x}=(x-x_k)/H_x$ and $\ul{t}=(t-t_k)/H_t$. Of course, many other choices for $w$ are possible too.

\begin{figure}[t]
\centering
\includegraphics[width=\linewidth,height=0.5\columnwidth]{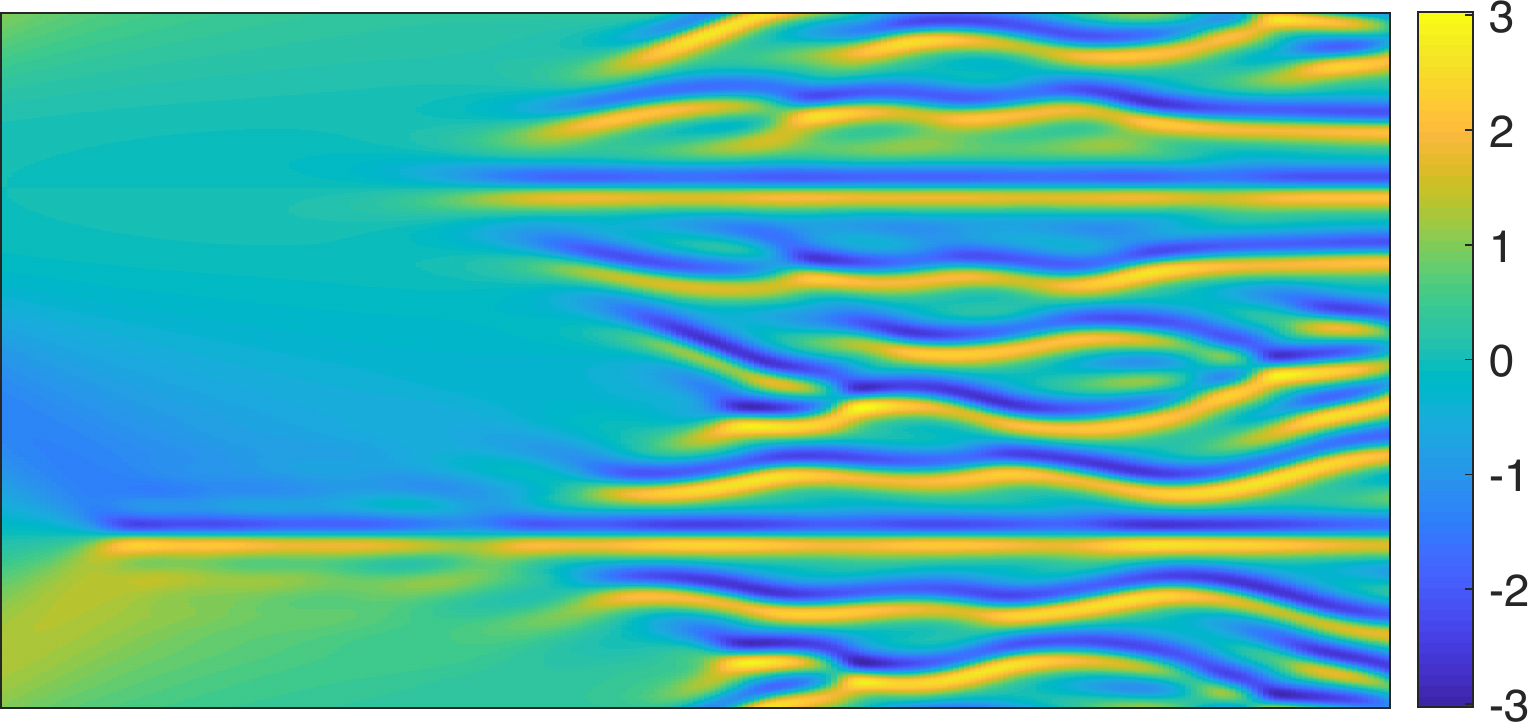}
\caption{Space-time plot of the solution to the Kuramoto-Sivashinsky equation. The $x$ axis is vertical and the $t$ axis is horizontal. The figures below use a similar colormap.}
\label{fig:KS_spacetime}
\end{figure}
 
The linear system \eqref{eq:symreg} can now be constructed by evaluating the integrals in \eqref{eq:KSlib} over a set of domains $\Omega_k$. 
To test our sparse regression approach, we generated surrogate data by solving the Kuramoto-Sivashinsky equation numerically.
To enable direct comparison with the results of Rudy {\it et al.} \cite{rudy2017}, we used the same integrator \cite{kassam2005} to compute the solution of \eqref{eq:KS} on a spatiotemporal domain of size $L_x=32\pi$ and $L_t=100$ using a grid with the same density $\Delta x=0.0982$ and $\Delta t=0.4$; the solution is shown in Fig. \ref{fig:KS_spacetime}.
Gaussian noise with standard deviation $\sigma$ was then added to $u$ at each grid point, after which the integrals in \eqref{eq:KSlib} were evaluated over integration domains with dimensions $H_x\approx24.5$, $H_t=20$.
The weight function used the exponents $p=4$ and $q=3$. 

\begin{figure}[t]
\centering
\includegraphics[width=0.9\columnwidth]{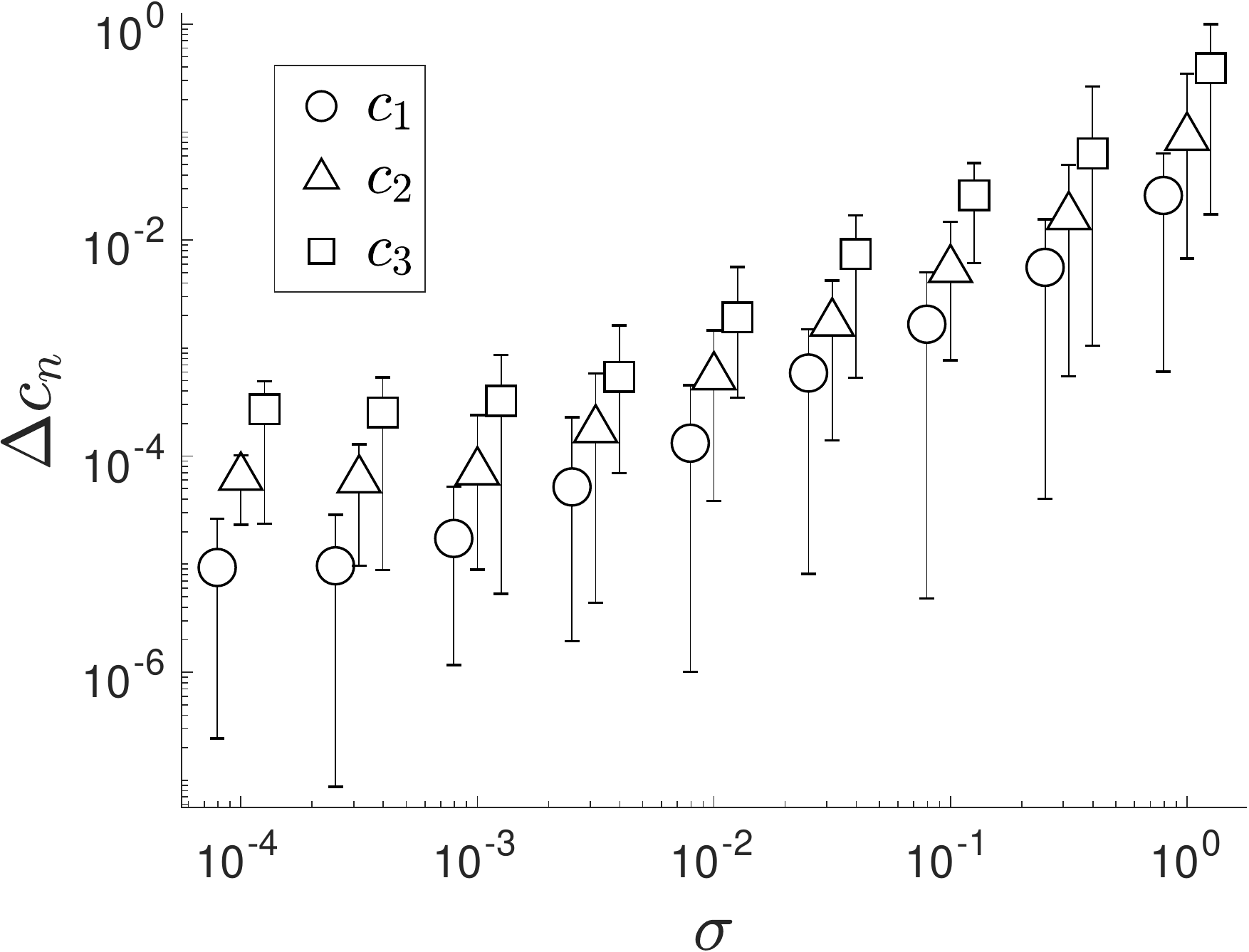}
\caption{The accuracy of parameter reconstruction for the Kuramoto-Sivashinsky equation as a function of the noise amplitude.}
\label{fig:KS_noise}
\end{figure}

The results for different noise levels are shown in Fig. \ref{fig:KS_noise}, with the accuracy of the model reconstruction quantified by the relative errors 
\begin{align}
\Delta c_n = \left|\frac{\tilde{c}_n-c_n}{c_n}\right|,
\label{eq:c_err}
\end{align}
where $c_n$ are the coefficients used to generate the numerical data and $\tilde{c}_n$ are the coefficients estimated from noisy data by via our sparse regression algorithm. 
Here and below, the symbols and the error bars show the mean values and the full range of the results, respectively, for the entire ensemble.
Note that the reconstruction remains essentially unaffected by noise, with error of about 1\% or below, until the noise level exceeds 10\%. 
This is a dramatic improvement compared to the original study \cite{rudy2017}, which yielded errors of over 50\% for all of the coefficients with just 1\% noise.

\section*{\label{sec:Kolmo} Latent variables}

To illustrate our approach applied to systems with latent variables, we next consider a flow in a thin layer of fluid driven by a steady but spatially nonuniform force ${\bf f}$.
The flow can be described using a generalization of the two-dimensional Navier-Stokes equation
\begin{align}
\partial_t{\bf u} = c_1({\bf u}\cdot\nabla){\bf u} + c_2\nabla^2{\bf u} + c_3{\bf u}-\nabla p + {\bf f},
\label{eq:NSE}
\end{align}
where ${\bf u}=\hat{x}u_x+\hat{y}u_y$ is the flow field, which is considered to be incompressible, $p$ is the pressure, and the constants $c_1$, $c_2$, and $c_3$ describe, respectively, the depth-averaged effects of inertia and viscosity in the horizontal and vertical direction \cite{suri_2014,tithof_2017}. In this example, both $p$ and ${\bf f}$ are assumed to be latent variables that cannot be measured.

To convert \eqref{eq:NSE} to weak form, we multiply it by a vector field ${\bf w}$ and integrate the result by parts over a (now three-dimensional) spatiotemporal domain $\Omega_k$ of size $2H_x\times2H_y\times2H_t$. Assuming again that the boundary terms vanish, for the linear terms we immediately find
\begin{align}
q^k_0  & = \int_{\Omega_k} {\bf w}\cdot\partial_t{\bf u}\,d\Omega
=-\int_{\Omega_k} {\bf u}\cdot\partial_t {\bf w} \,d\Omega, \nonumber\\
q^k_2 & = \int_{\Omega_k} {\bf w}\cdot\nabla^2{\bf u} \,d\Omega=
\int_{\Omega_k} {\bf u}\cdot\nabla^2{\bf w} \,d\Omega, \nonumber\\
q^k_3 & = \int_{\Omega_k} {\bf w}\cdot{\bf u}\,d\Omega. 
\label{eq:NSlib1}
\end{align}
The nonlinear term can be rewritten in a similar way using the incompressibility condition $\partial_iu_i=0$ (where summation over repeated indices is implied):
\begin{align}
q^k_1  & = \int_{\Omega_k} w_iu_j\partial_ju_i\,d\Omega
=-\int_{\Omega_k} u_i\partial_j(w_iu_j) \,d\Omega\nonumber\\
&=-\int_{\Omega_k} u_iu_j\partial_jw_i \,d\Omega
=-\int_{\Omega_k} {\bf u}\cdot ({\bf u}\cdot\nabla){\bf w} \,d\Omega.
\label{eq:NSlib2}
\end{align}
Finally, for the terms involving the latent variables, we find
\begin{align}
q^k_4 & = \int_{\Omega_k} {\bf w}\cdot\nabla p\,d\Omega
=-\int_{\Omega_k} p\nabla\cdot {\bf w} \,d\Omega, \nonumber\\
q^k_5 & = \int_{\Omega_k} {\bf w}\cdot{\bf f} \,d\Omega. 
\label{eq:NSlib3}
\end{align}

In order for the boundary terms to vanish on a rectangular domain $\Omega_k$ centered at $(x_k,y_k,t_k)$, we need to have ${\bf w}=0$ on $\partial\Omega$, as well as $\partial_x{\bf w}=0$ at $\ul{x}=\pm 1$ and $\partial_y{\bf w}=0$ at $\ul{y}=\pm 1$, where the underbar denotes rescaled variables $\ul{x}=(x-x_k)/H_x$, $\ul{y}=(y-y_k)/H_y$, and $\ul{t}=(t-t_k)/H_t$.
Next, the dependence on the pressure field and the steady forcing can be eliminated by additionally requiring that
\begin{align}
\nabla\cdot{\bf w}=0
\end{align}
and 
\begin{align}
\int_{-1}^1{\bf w}\,d\ul{t}=0.
\end{align}
All of the above conditions on ${\bf w}$ can be satisfied by setting ${\bf w}=\nabla\times (\psi\hat{z})=\hat{x}\partial_y\psi-\hat{y}\partial_x\psi$, where 
\begin{align}
\psi=\sin(\pi\ul{t})(\ul{x}^2-1)^p(\ul{y}^2-1)^p
\end{align} 
and $p\ge3$ (we used $p=3$ in this study). This yields $q^k_4=q^k_5=0$ and
\begin{align}
q^k_0  & = -\int_{\Omega_k} (u_x\partial_y - u_y\partial_x)\partial_t\psi\,d\Omega, \nonumber\\
q^k_1 & = \int_{\Omega_k} \left[ u_xu_y(\partial^2_y - \partial^2_x) + (u_x^2-u_y^2)\partial_{xy}\right]\psi\,d\Omega, \nonumber\\
q^k_2 & = \int_{\Omega_k} (u_x\partial_y - u_y\partial_x)\nabla^2\psi\,d\Omega, \nonumber\\
q^k_3 & = \int_{\Omega_k} (u_x\partial_y - u_y\partial_x)\psi\,d\Omega. 
\label{eq:NSlib}
\end{align} 
As in the case of the Kuramoto-Sivashinsky equation, the linear system \eqref{eq:symreg} can now be constructed by evaluating the integrals in \eqref{eq:NSlib} over a set of domains $\Omega_k$.
Note that this linear system involves neither the derivatives of the noisy observable data (components of the ${\bf u}$ field) nor the latent variables ($p$ and ${\bf f}$ fields).

\begin{figure}[t]
\subfigure[]{\includegraphics[width=0.48\columnwidth]{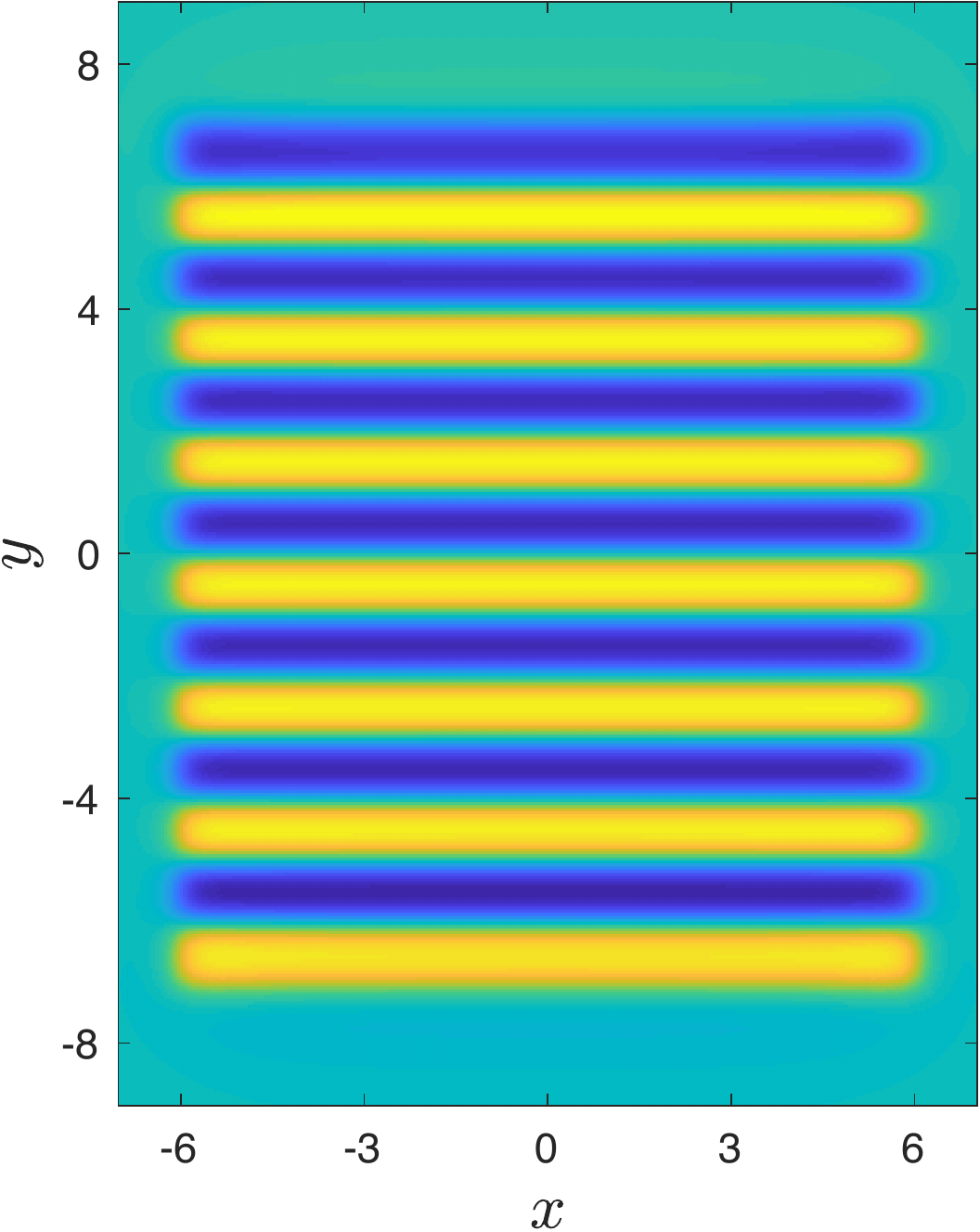}}\hspace{2mm}
\subfigure[]{\includegraphics[width=0.48\columnwidth]{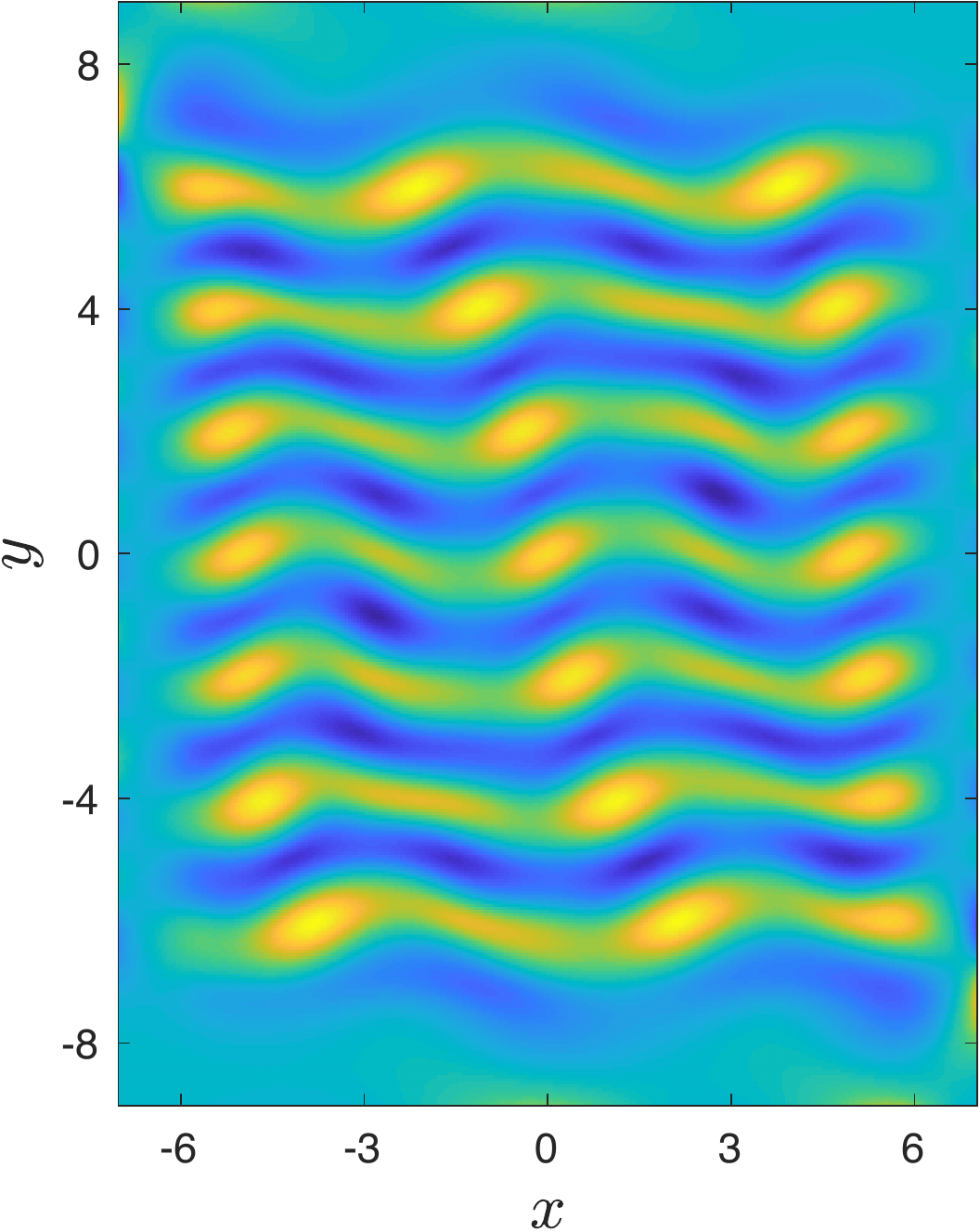}}
\caption{The Kolmogorov-like flow: (a) the forcing profile $f$ and (b) snapshot of the vorticity $\omega = \partial_xu_y-\partial_yu_x$. The $x$ axis is horizontal and the $y$ axis is vertical.}
\label{fig:Kolmo_forcing}
\end{figure}

To test our approach, we generated surrogate data ${\bf u}$ by solving \eqref{eq:NSE} with the parameters $c_1=-0.826$, $c_2=0.0487$, and $c_3=-0.157$, which correspond to the experimental setup of Kolmogorov-like flow described in Ref. \cite{tithof_2017}. In the experiment, the forcing field ${\bf f}=f(x,y)\hat{x}$ is produced by an array of long bar magnets with alternating polarity and width equal to unity in nondimensional units; correspondingly, $f(x,y)$ is approximately uniform in the $x$ direction and nearly periodic in the $y$ direction (cf. Fig. \ref{fig:Kolmo_forcing}(a)), with the ``period'' equal to 2 units. Forcing with amplitude $\max_{x,y}|f(x,y)|=1.0649$ generates a weakly turbulent flow (a representative snapshot is shown in Fig. \ref{fig:Kolmo_forcing}(b)), which was computed using the numerical integrator described in Ref. \cite{tithof_2017} on a domain of size $L_x=14$, $L_y=18$, $L_t\approx920$ and a computational grid with $\Delta x_c=\Delta y_c=0.025$ and $\Delta t_c\approx0.02$. 
The data was then subsampled on a coarser grid with spacing $\Delta x=\Delta y=0.1$ and $\Delta t=0.2302$, and Gaussian random noise with variance $\sigma$ was added to both components of the flow velocity ${\bf u}$. 
The integrals in \eqref{eq:NSlib} were evaluated over domains $\Omega_k$ of size $H_x=11.2$, $H_y=14.4$, and $H_t\approx34.5$.

\begin{figure}[t]
\centering
\includegraphics[width=0.9\columnwidth]{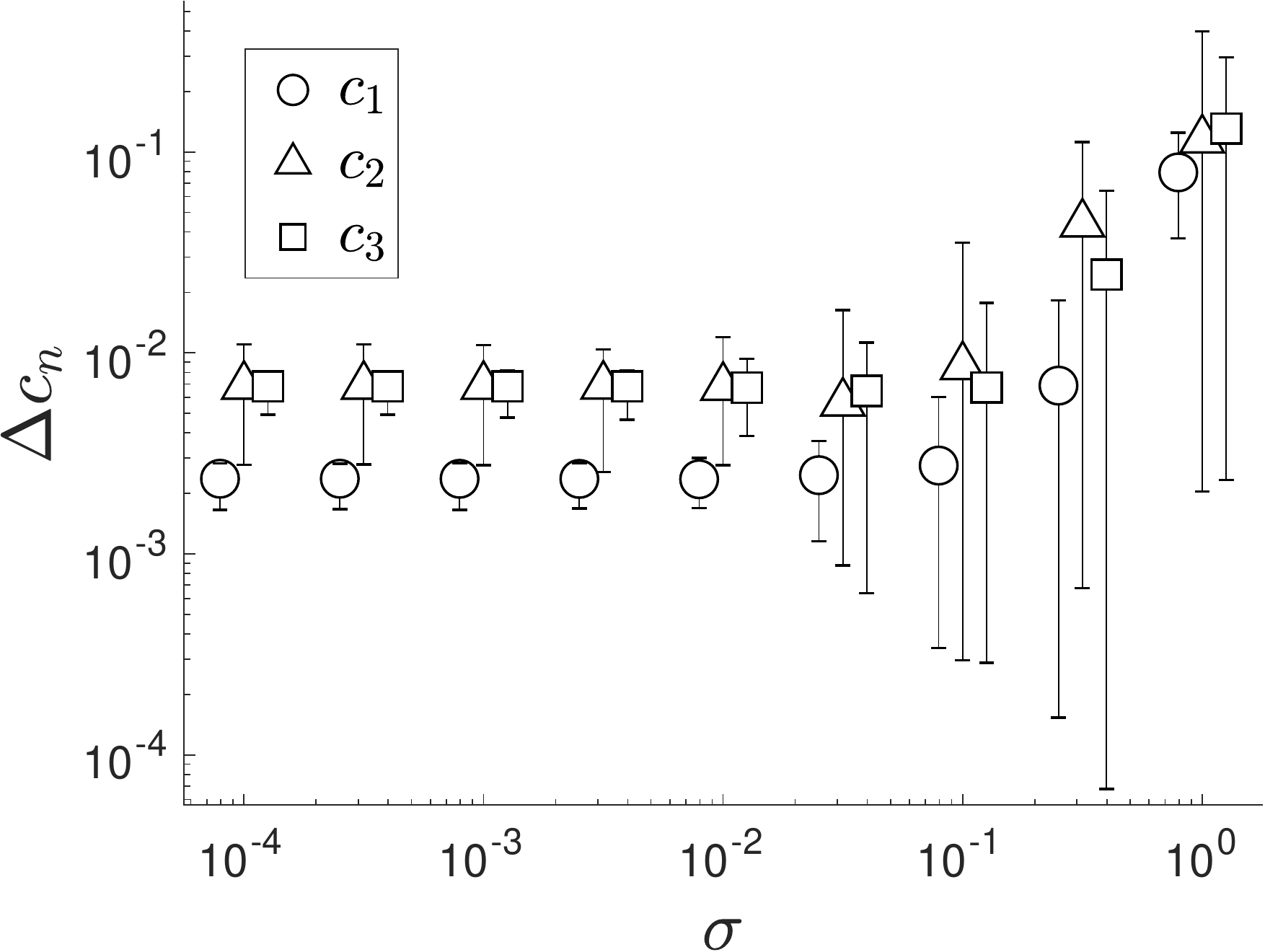}
\caption{The accuracy of parameter reconstruction for the 2D Kolmogorov-like flow model as a function of the noise amplitude.}
\label{fig:NS_noise}
\end{figure}

As Fig. \ref{fig:NS_noise} illustrates, our approach successfully reconstructs the reference PDE \eqref{eq:NSE}. 
Just like in the case of the Kuramoto-Sivashinsky equation, noise up to 10\% does not meaningfully affect the accuracy of model reconstruction, with the coefficients $c_1$, $c_2$, and $c_3$ estimated to within 1\% or better.
In fact, even with 100\% noise, the coefficients can still be estimated to within roughly 10\%.
For reference, experimental data \cite{tithof_2017} obtained using particle image velocimetry has roughly 3\% noise, at which level local sparse regression \cite{reinbold_2019} failed completely. 

\section*{\label{sec:RD} Sparse regression}

Finally, as an example of how the proposed approach could be used in the context of sparse regression, we consider the $\lambda-\omega$ reaction-diffusion system \cite{kopell_1973} in two spatial dimensions,
\begin{align}
\partial_tu &= D\nabla^2u+\lambda u-\omega v,\nonumber\\
\partial_tv &= D\nabla^2v+\omega u+\lambda v,
\label{eq:RDE}
\end{align}
where $\omega=-\beta (u^2+v^2)$, $\lambda=1-u^2-v^2$, and $\beta=1$ and $D=0.1$ 
are constants. 
This system can be cast in the form of Eq.~\eqref{eq:nl2} by defining a vector ${\bf u}=[u,v]$. To test our approach, we applied sparse regression to a generalization of \eqref{eq:RDE}, where the reaction terms are given by polynomials in $u$ and $v$ up to third order. In total, the generalized model involves a total of 20 different terms (two diffusion terms and 18 polynomial terms). Correspondingly, 20 unknown coefficients need to be determined. 

The sparse regression problem for the $\lambda-\omega$ system can be block-diagonalized by using a weight function ${\bf w}=[w,0]$ (or ${\bf w}=[0,w]$) to reconstruct the first (or second) equation in \eqref{eq:RDE}, yielding two independent linear systems \eqref{eq:symreg} with 10 library terms each. 
The integration domains $\Omega_k$ are three-dimensional as in the previous example.
The integrals involving terms such as $u^\alpha v^\beta$ do not require integration by parts.
The two integrals involving the Laplacian terms are integrated by parts twice to get rid of derivatives on $u$ and $v$, e.g.,
\begin{align}
q_1^k = \int_{\Omega_k} w \nabla^2u\,d\Omega=\int_{\Omega_k} u \nabla^2w\,d\Omega.
\end{align}
In both cases, the corresponding boundary terms vanish if we choose 
\begin{equation}
w = (\ul{x}^2-1)^p(\ul{y}^2-1)^p(\ul{t}^2-1)^q,
\end{equation}
where $p\ge2$ and $q\ge1$ (we chose $p=2$ and $q=1$). 

\begin{figure}[t]
\centering
\subfigure[]{\includegraphics[width=0.48\columnwidth]{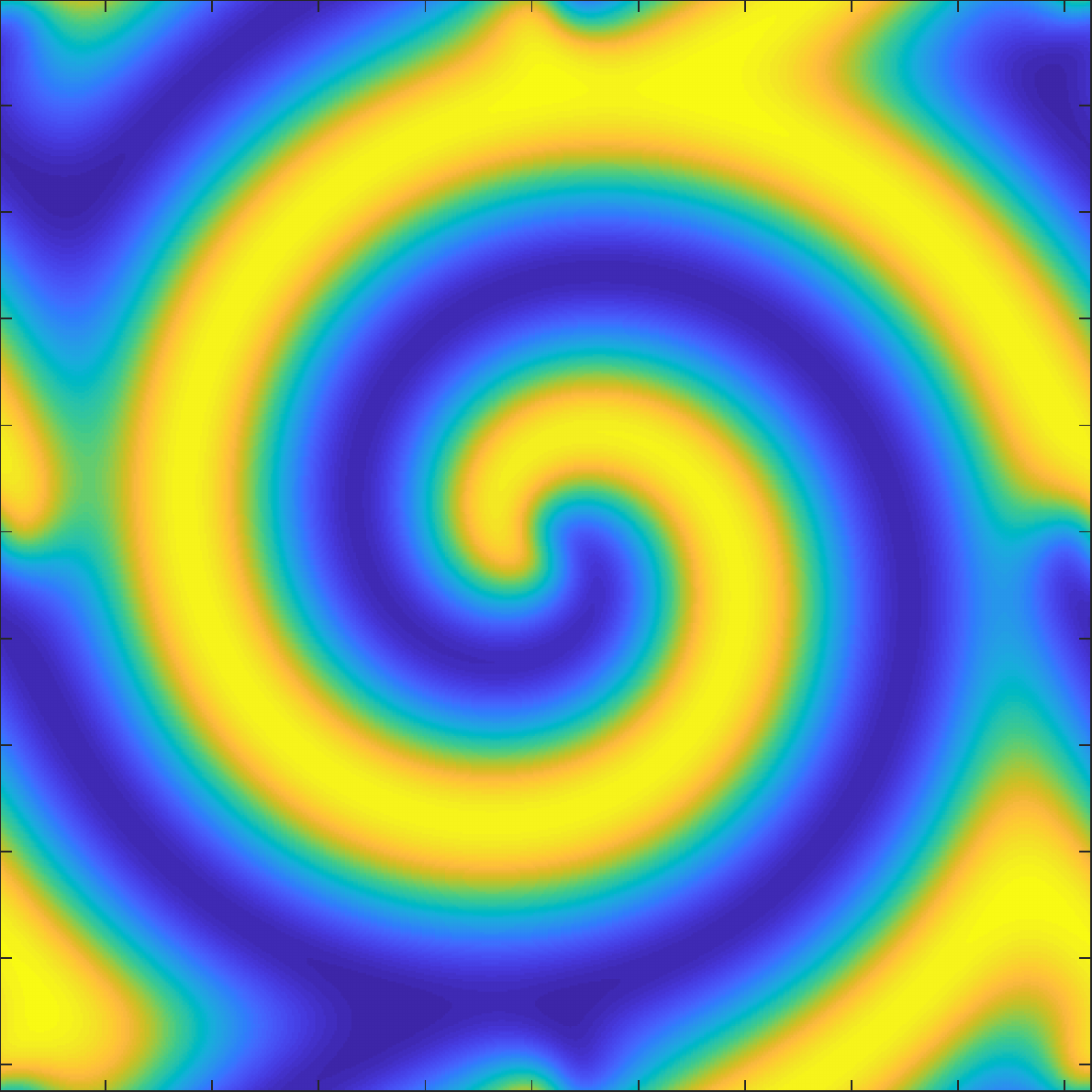}}\hspace{2mm}
\subfigure[]{\includegraphics[width=0.48\columnwidth]{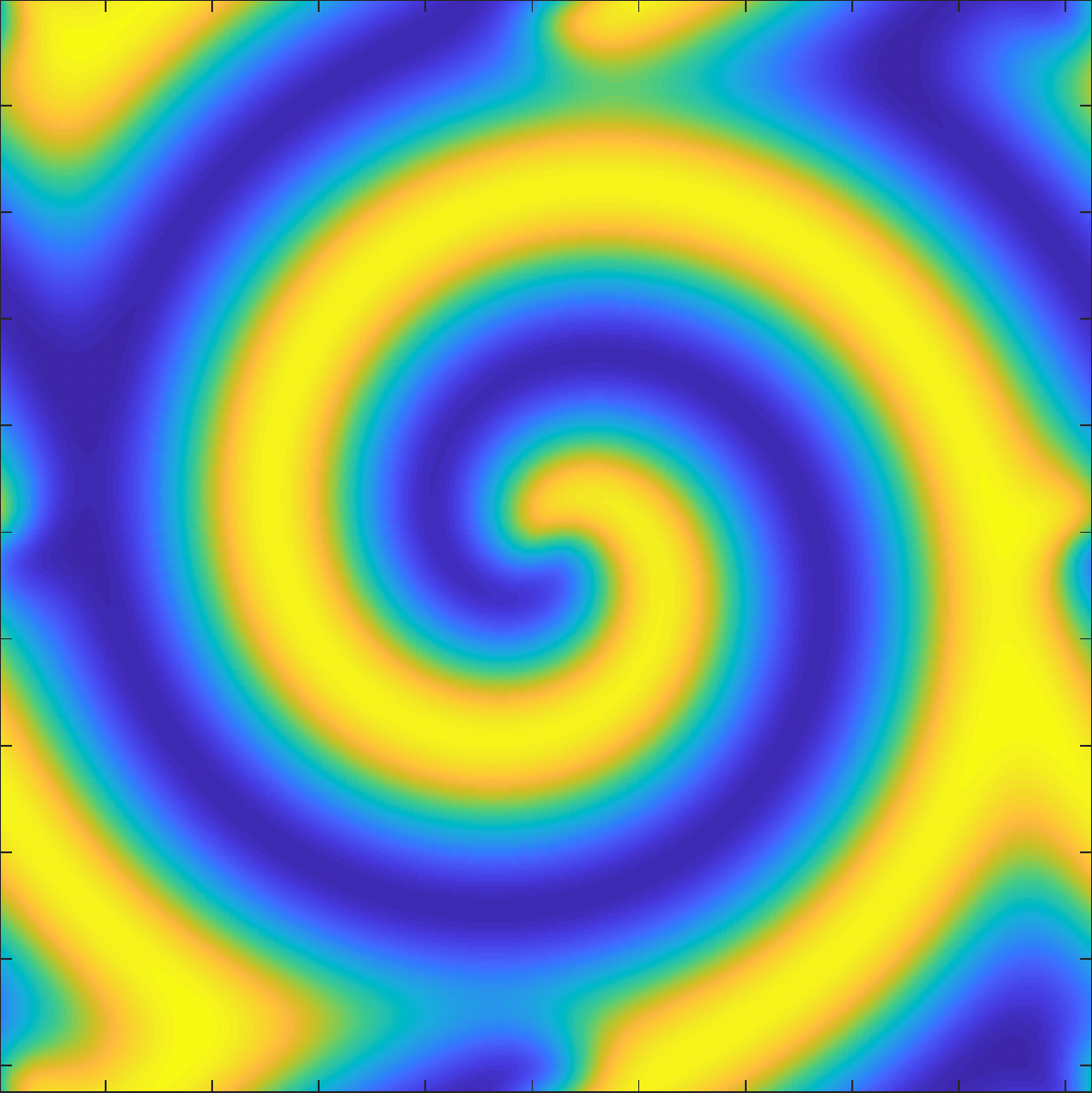}}
\caption{A typical snapshot of the fields (a) $u$ and (b) $v$ for the $\lambda$-$\omega$ reaction diffusion system.  The $x$ axis is horizontal and the $y$ axis is vertical.}
\label{fig:RD}
\end{figure}

The surrogate data was obtained by computing the solution of \eqref{eq:RDE} using the integrator employed in Ref. \cite{rudy2017}; a typical snapshot is shown in Fig. \ref{fig:RD}. 
The computational domain of size $L_x=20$, $L_y=20$, $L_t=10$ was discretized using a grid with spacing $\Delta x=\Delta y=0.0391$ and $\Delta t=0.05$, and Gaussian random noise with standard deviation $\sigma$ was added to both $u$ and $v$ at each grid point. 
The dimensions of the integration domains $\Omega_k$ were chosen as $H_x=H_y\approx1$ and $H_t=1.25$. 

\begin{figure}[t]
\centering
\includegraphics[width=0.9\columnwidth]{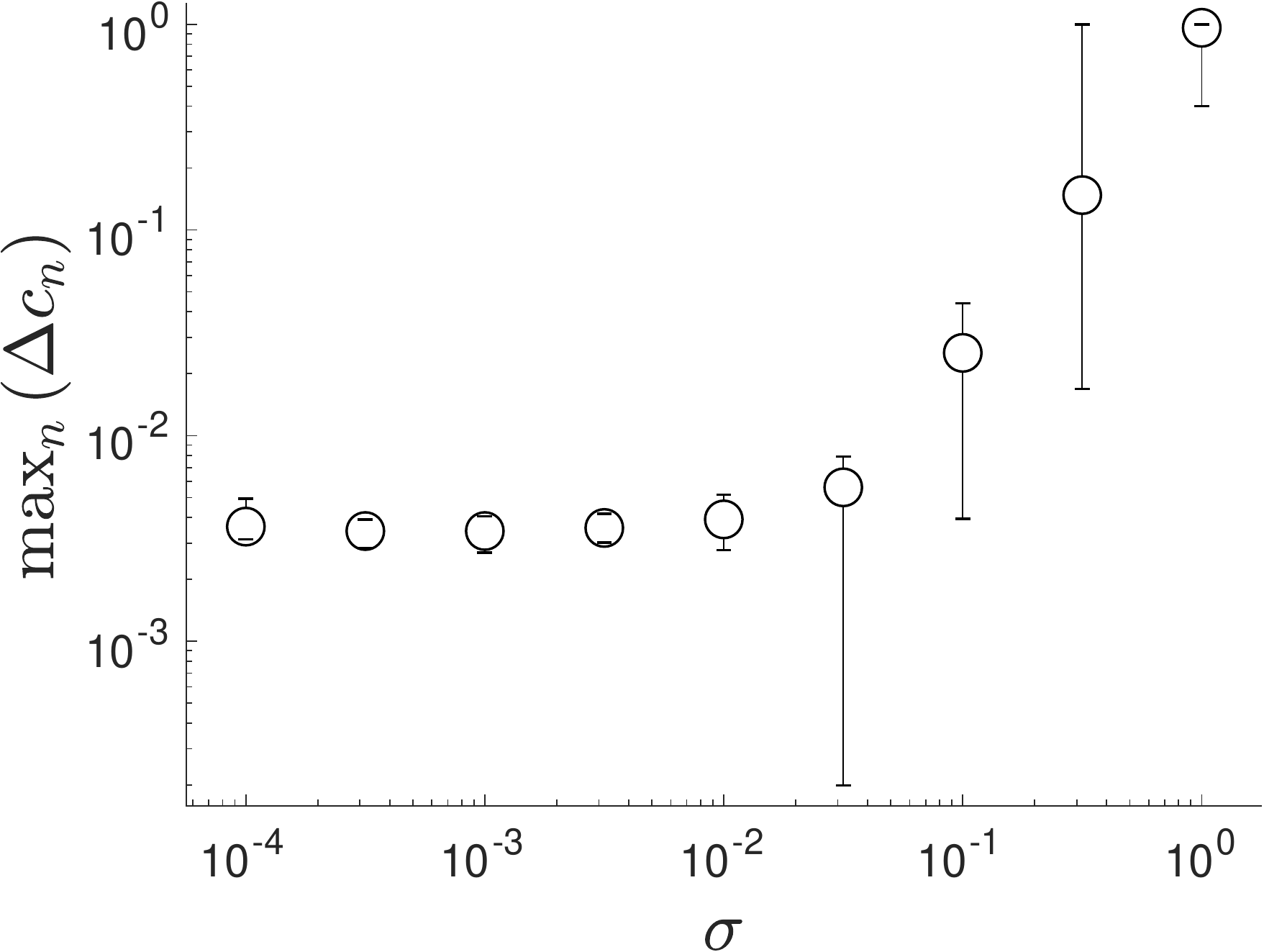}
\caption{The accuracy of parameter reconstruction for the $\lambda-\omega$ reaction-diffusion system as a function of the noise amplitude. Shown is the largest error, which corresponds to one of the diffusion coefficients. 
}
\label{fig:RD_noise}
\end{figure}

The results of sparse regression are shown in Fig. \ref{fig:RD_noise}. We find that, for noise levels of up to 5\%, the model was reconstructed correctly (with no spurious or missing terms) for each distribution of $\Omega_k$ in our ensemble, with all parameters estimated to an accuracy of better than 1\%. With 10\% noise, the model is identified correctly in about 95\% of cases, and at 30\% noise, the model is identified correctly in about 20\% of cases, with the remaining cases featuring spurious terms (linear in $u$ and $v$) that are not present in the $\lambda-\omega$ model. For reference, sparse regression based on local evaluation of derivatives \cite{rudy2017} failed to correctly identify this model, generating spurious terms in the presence of as little as 1\% noise.

It should be noted that using ensemble sparse regression makes it easy to detect the presence of spurious (missing) terms and eliminate (add) them while still preserving the accuracy with which all of the correct terms are estimated (in our case, about 3\% for the worst case offenders with 10\% noise). It is also worth pointing out that, unlike the standard approach \cite{rudy2017}, weak formulation requires no intermediate noise reduction.

\section*{Discussion}

The examples presented here illustrate the power of the weak formulation of sparse regression applied to noisy and/or incomplete data. For instance, high-order PDEs such as the Kuramoto-Sivashinsky equation simply cannot be reconstructed with meaningful accuracy from data with realistic levels of noise using the original (differential) form of the model. The main culprit is the term in the model involving a fourth-order derivative, which is extremely sensitive even to minute amounts of noise. The weak formulation involves integrals of the data rather than derivatives, which makes it much more robust with respect to noise. While the weak formulation may not eliminate {\it all} of the derivatives in some models (e.g., in nonlinear terms), it can reduce the order of the derivatives that remain, which is extremely beneficial when noisy data is involved. 

We have also demonstrated that the weak formulation of sparse regression can be applied successfully to models with latent variables, as in the example of the fluid flow in a thin layer, where neither the pressure field nor the forcing field are accessible. Needless to say, the weak formulation by itself simply eliminates rather than reconstructs the terms that involve the latent variables. One needs to impose additional physical constraints \cite{reinbold_2019} to determine their functional form. Nonetheless, the approach presented here has substantial advantages compared to the method described in Ref. \cite{reinbold_2019}, which involves taking additional spatial and/or temporal derivatives of the model equation to eliminate the latent variables. As discussed previously, the higher the order of the derivatives, the more sensitive the sparse regression is to noise. As a result, the model \eqref{eq:NSE} could only be reconstructed with acceptable accuracy in that study for noise levels of 0.01\% or less. The present approach gives better accuracy for data with as much as 30\% noise!

In conclusion, let us point out that we have made no attempt to optimize our approach here. Several options are available to make it even more robust and accurate \cite{gurevich_2019}. As an example, the size of the integration domains $\Omega_k$ could be varied relative to the size of the spatiotemporal domain on which the data are available. Furthermore, we have only used a single weight function, while in principle one could also use a set of different weight functions ${\bf w}_j$. Additionally, the shape of the weight functions could be optimized to improve the accuracy even compared to the already impressive results presented here. For instance, simply increasing the powers $p$ and $q$ beyond the minimal possible values (determined, respectively, by the highest order of the spatial and temporal derivatives in the model) can reduce the error in estimating the coefficients of the model by orders of magnitude.
In contrast, we found the details of the sparse regression procedure itself to have a relatively minor impact on the results.

\acknowledgements{This material is based upon work supported by NSF under Grant No. CMMI-1725587 and ARO under grant No. W911NF-15-10471.  DRG gratefully acknowledges the support of the Letson Undergraduate Research Scholarship.}


\begin{thebibliography}{19}%
\makeatletter
\providecommand \@ifxundefined [1]{%
 \@ifx{#1\undefined}
}%
\providecommand \@ifnum [1]{%
 \ifnum #1\expandafter \@firstoftwo
 \else \expandafter \@secondoftwo
 \fi
}%
\providecommand \@ifx [1]{%
 \ifx #1\expandafter \@firstoftwo
 \else \expandafter \@secondoftwo
 \fi
}%
\providecommand \natexlab [1]{#1}%
\providecommand \enquote  [1]{``#1''}%
\providecommand \bibnamefont  [1]{#1}%
\providecommand \bibfnamefont [1]{#1}%
\providecommand \citenamefont [1]{#1}%
\providecommand \href@noop [0]{\@secondoftwo}%
\providecommand \href [0]{\begingroup \@sanitize@url \@href}%
\providecommand \@href[1]{\@@startlink{#1}\@@href}%
\providecommand \@@href[1]{\endgroup#1\@@endlink}%
\providecommand \@sanitize@url [0]{\catcode `\\12\catcode `\$12\catcode
  `\&12\catcode `\#12\catcode `\^12\catcode `\_12\catcode `\%12\relax}%
\providecommand \@@startlink[1]{}%
\providecommand \@@endlink[0]{}%
\providecommand \url  [0]{\begingroup\@sanitize@url \@url }%
\providecommand \@url [1]{\endgroup\@href {#1}{\urlprefix }}%
\providecommand \urlprefix  [0]{URL }%
\providecommand \Eprint [0]{\href }%
\providecommand \doibase [0]{http://dx.doi.org/}%
\providecommand \selectlanguage [0]{\@gobble}%
\providecommand \bibinfo  [0]{\@secondoftwo}%
\providecommand \bibfield  [0]{\@secondoftwo}%
\providecommand \translation [1]{[#1]}%
\providecommand \BibitemOpen [0]{}%
\providecommand \bibitemStop [0]{}%
\providecommand \bibitemNoStop [0]{.\EOS\space}%
\providecommand \EOS [0]{\spacefactor3000\relax}%
\providecommand \BibitemShut  [1]{\csname bibitem#1\endcsname}%
\let\auto@bib@innerbib\@empty
\bibitem [{\citenamefont {Crutchfield}\ and\ \citenamefont
  {McNamara}(1987)}]{crutchfield1987}%
  \BibitemOpen
  \bibfield  {author} {\bibinfo {author} {\bibfnamefont {J.~P.}\ \bibnamefont
  {Crutchfield}}\ and\ \bibinfo {author} {\bibfnamefont {B.~S.}\ \bibnamefont
  {McNamara}},\ }\href@noop {} {\bibfield  {journal} {\bibinfo  {journal}
  {Complex systems}\ }\textbf {\bibinfo {volume} {1}},\ \bibinfo {pages} {417}
  (\bibinfo {year} {1987})}\BibitemShut {NoStop}%
\bibitem [{\citenamefont {Bongard}\ and\ \citenamefont
  {Lipson}(2007)}]{bongard2007}%
  \BibitemOpen
  \bibfield  {author} {\bibinfo {author} {\bibfnamefont {J.}~\bibnamefont
  {Bongard}}\ and\ \bibinfo {author} {\bibfnamefont {H.}~\bibnamefont
  {Lipson}},\ }\href@noop {} {\bibfield  {journal} {\bibinfo  {journal}
  {Proceedings of the National Academy of Sciences}\ }\textbf {\bibinfo
  {volume} {104}},\ \bibinfo {pages} {9943} (\bibinfo {year}
  {2007})}\BibitemShut {NoStop}%
\bibitem [{\citenamefont {Yao}\ and\ \citenamefont {Bollt}(2007)}]{yao2007}%
  \BibitemOpen
  \bibfield  {author} {\bibinfo {author} {\bibfnamefont {C.}~\bibnamefont
  {Yao}}\ and\ \bibinfo {author} {\bibfnamefont {E.~M.}\ \bibnamefont
  {Bollt}},\ }\href@noop {} {\bibfield  {journal} {\bibinfo  {journal} {Physica
  D: Nonlinear Phenomena}\ }\textbf {\bibinfo {volume} {227}},\ \bibinfo
  {pages} {78} (\bibinfo {year} {2007})}\BibitemShut {NoStop}%
\bibitem [{\citenamefont {Chou}\ and\ \citenamefont {Voit}(2009)}]{chou2009}%
  \BibitemOpen
  \bibfield  {author} {\bibinfo {author} {\bibfnamefont {I.-C.}\ \bibnamefont
  {Chou}}\ and\ \bibinfo {author} {\bibfnamefont {E.~O.}\ \bibnamefont
  {Voit}},\ }\href@noop {} {\bibfield  {journal} {\bibinfo  {journal}
  {Mathematical biosciences}\ }\textbf {\bibinfo {volume} {219}},\ \bibinfo
  {pages} {57} (\bibinfo {year} {2009})}\BibitemShut {NoStop}%
\bibitem [{\citenamefont {Brunton}\ \emph {et~al.}(2016)\citenamefont
  {Brunton}, \citenamefont {Proctor},\ and\ \citenamefont
  {Kutz}}]{brunton2016}%
  \BibitemOpen
  \bibfield  {author} {\bibinfo {author} {\bibfnamefont {S.~L.}\ \bibnamefont
  {Brunton}}, \bibinfo {author} {\bibfnamefont {J.~L.}\ \bibnamefont
  {Proctor}}, \ and\ \bibinfo {author} {\bibfnamefont {J.~N.}\ \bibnamefont
  {Kutz}},\ }\href@noop {} {\bibfield  {journal} {\bibinfo  {journal}
  {Proceedings of the National Academy of Sciences}\ }\textbf {\bibinfo
  {volume} {113}},\ \bibinfo {pages} {3932} (\bibinfo {year}
  {2016})}\BibitemShut {NoStop}%
\bibitem [{\citenamefont {Xu}\ and\ \citenamefont
  {Khanmohamadi}(2008)}]{xu_2008}%
  \BibitemOpen
  \bibfield  {author} {\bibinfo {author} {\bibfnamefont {D.}~\bibnamefont
  {Xu}}\ and\ \bibinfo {author} {\bibfnamefont {O.}~\bibnamefont
  {Khanmohamadi}},\ }\href@noop {} {\bibfield  {journal} {\bibinfo  {journal}
  {Chaos}\ }\textbf {\bibinfo {volume} {18}},\ \bibinfo {pages} {043122}
  (\bibinfo {year} {2008})}\BibitemShut {NoStop}%
\bibitem [{\citenamefont {Rudy}\ \emph {et~al.}(2017)\citenamefont {Rudy},
  \citenamefont {Brunton}, \citenamefont {Proctor},\ and\ \citenamefont
  {Kutz}}]{rudy2017}%
  \BibitemOpen
  \bibfield  {author} {\bibinfo {author} {\bibfnamefont {S.~H.}\ \bibnamefont
  {Rudy}}, \bibinfo {author} {\bibfnamefont {S.~L.}\ \bibnamefont {Brunton}},
  \bibinfo {author} {\bibfnamefont {J.~L.}\ \bibnamefont {Proctor}}, \ and\
  \bibinfo {author} {\bibfnamefont {J.~N.}\ \bibnamefont {Kutz}},\ }\href@noop
  {} {\bibfield  {journal} {\bibinfo  {journal} {Science Advances}\ }\textbf
  {\bibinfo {volume} {3}},\ \bibinfo {pages} {e1602614} (\bibinfo {year}
  {2017})}\BibitemShut {NoStop}%
\bibitem [{\citenamefont {Li}\ \emph {et~al.}(2019)\citenamefont {Li},
  \citenamefont {Li}, \citenamefont {Yue}, \citenamefont {Tang}, \citenamefont
  {Voss}, \citenamefont {Kurths},\ and\ \citenamefont {Yuan}}]{li_2019}%
  \BibitemOpen
  \bibfield  {author} {\bibinfo {author} {\bibfnamefont {X.}~\bibnamefont
  {Li}}, \bibinfo {author} {\bibfnamefont {L.}~\bibnamefont {Li}}, \bibinfo
  {author} {\bibfnamefont {Z.}~\bibnamefont {Yue}}, \bibinfo {author}
  {\bibfnamefont {X.}~\bibnamefont {Tang}}, \bibinfo {author} {\bibfnamefont
  {H.~U.}\ \bibnamefont {Voss}}, \bibinfo {author} {\bibfnamefont
  {J.}~\bibnamefont {Kurths}}, \ and\ \bibinfo {author} {\bibfnamefont
  {Y.}~\bibnamefont {Yuan}},\ }\href@noop {} {\bibfield  {journal} {\bibinfo
  {journal} {Chaos}\ }\textbf {\bibinfo {volume} {29}},\ \bibinfo {pages}
  {043130} (\bibinfo {year} {2019})}\BibitemShut {NoStop}%
\bibitem [{\citenamefont {Reinbold}\ and\ \citenamefont
  {Grigoriev}(2019)}]{reinbold_2019}%
  \BibitemOpen
  \bibfield  {author} {\bibinfo {author} {\bibfnamefont {P.~A.~K.}\
  \bibnamefont {Reinbold}}\ and\ \bibinfo {author} {\bibfnamefont {R.~O.}\
  \bibnamefont {Grigoriev}},\ }\href {\doibase 10.1103/PhysRevE.100.022219}
  {\bibfield  {journal} {\bibinfo  {journal} {Phys. Rev. E}\ }\textbf {\bibinfo
  {volume} {100}},\ \bibinfo {pages} {022219} (\bibinfo {year}
  {2019})}\BibitemShut {NoStop}%
\bibitem [{\citenamefont {Raissi}\ \emph {et~al.}(2018)\citenamefont {Raissi},
  \citenamefont {Perdikaris},\ and\ \citenamefont
  {Karniadakis}}]{raissi2018siam}%
  \BibitemOpen
  \bibfield  {author} {\bibinfo {author} {\bibfnamefont {M.}~\bibnamefont
  {Raissi}}, \bibinfo {author} {\bibfnamefont {P.}~\bibnamefont {Perdikaris}},
  \ and\ \bibinfo {author} {\bibfnamefont {G.~E.}\ \bibnamefont
  {Karniadakis}},\ }\href@noop {} {\bibfield  {journal} {\bibinfo  {journal}
  {SIAM Journal on Scientific Computing}\ }\textbf {\bibinfo {volume} {40}},\
  \bibinfo {pages} {A172} (\bibinfo {year} {2018})}\BibitemShut {NoStop}%
\bibitem [{\citenamefont {Raissi}\ and\ \citenamefont
  {Karniadakis}(2018)}]{raissi_2018}%
  \BibitemOpen
  \bibfield  {author} {\bibinfo {author} {\bibfnamefont {M.}~\bibnamefont
  {Raissi}}\ and\ \bibinfo {author} {\bibfnamefont {G.~E.}\ \bibnamefont
  {Karniadakis}},\ }\href@noop {} {\bibfield  {journal} {\bibinfo  {journal}
  {Journal of Computational Physics}\ }\textbf {\bibinfo {volume} {357}},\
  \bibinfo {pages} {125} (\bibinfo {year} {2018})}\BibitemShut {NoStop}%
\bibitem [{\citenamefont {Sivashinsky}(1977)}]{sivashinsky_1977}%
  \BibitemOpen
  \bibfield  {author} {\bibinfo {author} {\bibfnamefont {G.}~\bibnamefont
  {Sivashinsky}},\ }\href@noop {} {\bibfield  {journal} {\bibinfo  {journal}
  {Acta astronautica}\ }\textbf {\bibinfo {volume} {4}},\ \bibinfo {pages}
  {1177} (\bibinfo {year} {1977})}\BibitemShut {NoStop}%
\bibitem [{\citenamefont {Kuramoto}\ and\ \citenamefont
  {Tsuzuki}(1976)}]{kuramoto_1976}%
  \BibitemOpen
  \bibfield  {author} {\bibinfo {author} {\bibfnamefont {Y.}~\bibnamefont
  {Kuramoto}}\ and\ \bibinfo {author} {\bibfnamefont {T.}~\bibnamefont
  {Tsuzuki}},\ }\href@noop {} {\bibfield  {journal} {\bibinfo  {journal}
  {Progress of theoretical physics}\ }\textbf {\bibinfo {volume} {55}},\
  \bibinfo {pages} {356} (\bibinfo {year} {1976})}\BibitemShut {NoStop}%
\bibitem [{\citenamefont {Sivashinsky}\ and\ \citenamefont
  {Michelson}(1980)}]{sivashinsky_1980flow}%
  \BibitemOpen
  \bibfield  {author} {\bibinfo {author} {\bibfnamefont {G.~I.}\ \bibnamefont
  {Sivashinsky}}\ and\ \bibinfo {author} {\bibfnamefont {D.}~\bibnamefont
  {Michelson}},\ }\href@noop {} {\bibfield  {journal} {\bibinfo  {journal}
  {Progress of theoretical physics}\ }\textbf {\bibinfo {volume} {63}},\
  \bibinfo {pages} {2112} (\bibinfo {year} {1980})}\BibitemShut {NoStop}%
\bibitem [{\citenamefont {Kassam}\ and\ \citenamefont
  {Trefethen}(2005)}]{kassam2005}%
  \BibitemOpen
  \bibfield  {author} {\bibinfo {author} {\bibfnamefont {A.-K.}\ \bibnamefont
  {Kassam}}\ and\ \bibinfo {author} {\bibfnamefont {L.~N.}\ \bibnamefont
  {Trefethen}},\ }\href@noop {} {\bibfield  {journal} {\bibinfo  {journal}
  {SIAM Journal on Scientific Computing}\ }\textbf {\bibinfo {volume} {26}},\
  \bibinfo {pages} {1214} (\bibinfo {year} {2005})}\BibitemShut {NoStop}%
\bibitem [{\citenamefont {Suri}\ \emph {et~al.}(2014)\citenamefont {Suri},
  \citenamefont {Tithof}, \citenamefont {Mitchell}, \citenamefont {Grigoriev},\
  and\ \citenamefont {Schatz}}]{suri_2014}%
  \BibitemOpen
  \bibfield  {author} {\bibinfo {author} {\bibfnamefont {B.}~\bibnamefont
  {Suri}}, \bibinfo {author} {\bibfnamefont {J.}~\bibnamefont {Tithof}},
  \bibinfo {author} {\bibfnamefont {R.}~\bibnamefont {Mitchell}}, \bibinfo
  {author} {\bibfnamefont {R.~O.}\ \bibnamefont {Grigoriev}}, \ and\ \bibinfo
  {author} {\bibfnamefont {M.~F.}\ \bibnamefont {Schatz}},\ }\href {\doibase
  http://dx.doi.org/10.1063/1.4873417} {\bibfield  {journal} {\bibinfo
  {journal} {Phys. Fluids}\ }\textbf {\bibinfo {volume} {26}},\ \bibinfo {eid}
  {053601} (\bibinfo {year} {2014})}\BibitemShut {NoStop}%
\bibitem [{\citenamefont {Tithof}\ \emph {et~al.}(2017)\citenamefont {Tithof},
  \citenamefont {Suri}, \citenamefont {Pallantla}, \citenamefont {Grigoriev},\
  and\ \citenamefont {Schatz}}]{tithof_2017}%
  \BibitemOpen
  \bibfield  {author} {\bibinfo {author} {\bibfnamefont {J.}~\bibnamefont
  {Tithof}}, \bibinfo {author} {\bibfnamefont {B.}~\bibnamefont {Suri}},
  \bibinfo {author} {\bibfnamefont {R.~K.}\ \bibnamefont {Pallantla}}, \bibinfo
  {author} {\bibfnamefont {R.~O.}\ \bibnamefont {Grigoriev}}, \ and\ \bibinfo
  {author} {\bibfnamefont {M.~F.}\ \bibnamefont {Schatz}},\ }\href@noop {}
  {\bibfield  {journal} {\bibinfo  {journal} {Journal of Fluid Mechanics.}\
  }\textbf {\bibinfo {volume} {828}},\ \bibinfo {pages} {837} (\bibinfo {year}
  {2017})}\BibitemShut {NoStop}%
\bibitem [{\citenamefont {Kopell}\ and\ \citenamefont
  {Howard}(1973)}]{kopell_1973}%
  \BibitemOpen
  \bibfield  {author} {\bibinfo {author} {\bibfnamefont {N.}~\bibnamefont
  {Kopell}}\ and\ \bibinfo {author} {\bibfnamefont {L.~N.}\ \bibnamefont
  {Howard}},\ }\href@noop {} {\bibfield  {journal} {\bibinfo  {journal}
  {Studies in Applied Mathematics}\ }\textbf {\bibinfo {volume} {52}},\
  \bibinfo {pages} {291} (\bibinfo {year} {1973})}\BibitemShut {NoStop}%
\bibitem [{\citenamefont {Gurevich}\ \emph {et~al.}(2019)\citenamefont
  {Gurevich}, \citenamefont {Reinbold},\ and\ \citenamefont
  {Grigoriev}}]{gurevich_2019}%
  \BibitemOpen
  \bibfield  {author} {\bibinfo {author} {\bibfnamefont {D.~R.}\ \bibnamefont
  {Gurevich}}, \bibinfo {author} {\bibfnamefont {P.~A.~K.}\ \bibnamefont
  {Reinbold}}, \ and\ \bibinfo {author} {\bibfnamefont {R.~O.}\ \bibnamefont
  {Grigoriev}},\ }\href@noop {} {\bibfield  {journal} {\bibinfo  {journal}
  {Chaos}\ }\textbf {\bibinfo {volume} {29}},\ \bibinfo {pages} {103113}
  (\bibinfo {year} {2019})}\BibitemShut {NoStop}%
\end{thebibliography}
%

\end{document}